\documentclass[12pt]{amsart}
\usepackage{amssymb}
\usepackage{amsfonts}
\usepackage{latexsym}
\usepackage[all]{xy}
\usepackage{amscd}
\usepackage{comment}

\newtheorem{theorem}{Theorem}[section]
\newtheorem{lemma}[theorem]{Lemma}

\newtheorem{corollary}[theorem]{Corollary}
\theoremstyle{definition}
\newtheorem{definition}[theorem]{Definition}
\newtheorem{example}[theorem]{Example}

\newtheorem{remark}[theorem]{Remark}
\newtheorem{note}[theorem]{Note}

\newcommand{\C}{\mathcal{C}}
\newcommand{\D}{\mathcal{D}}
\newcommand{\V}{\mathcal{V}}
\newcommand{\M}{\mathcal{M}}

\renewcommand{\S}{\mathcal{S}}
\newcommand{\id}{\text{id}}
\newcommand{\morita}{\thickapprox}

\begin{document}

\title{Categorical Morita Equivalence For Group-Theoretical Categories}

\author{Deepak Naidu}
\address{Department of Mathematics and Statistics,
University of New Hampshire, Durham, NH 03824, USA}
\email{dnaidu@unh.edu}


\begin{abstract}
A finite tensor category is called {\em pointed} if all its simple objects are invertible.
We find necessary and sufficient conditions for two pointed semisimple categories 
to be dual to each other with respect to a module category. 
Whenever the dual of a pointed semisimple category 
with respect to a module category is pointed, we give explicit formulas for 
the Grothendieck ring and for the associator of the dual. This leads to the definition
of {\em categorical Morita equivalence} on the set of all finite groups and
on the set of all pairs $(G, \, \omega)$, where $G$ is a finite group and $\omega \in 
H^3(G, \, k^\times)$. A group-theoretical and cohomological interpretation of this 
relation is given. A series of concrete examples of pairs of groups that are 
categorically Morita equivalent but have non-isomorphic Grothendieck rings are given. 
In particular, the representation categories of 
the Drinfeld doubles of the groups in each example 
are equivalent as braided tensor categories and hence these groups define 
the same modular data.
\end{abstract}

\maketitle


\begin{section}
{Introduction}

Throughout this paper we work over an algebraically closed field $k$ 
of characteristic $0$. As pointed out to us by the referee, most results
can be extended to the positive characteristic setting using the theory
of finite tensor categories.
A right module category over a tensor category $\C$ is a semisimple 
abelian category $\M$ together with a functor $\M \times \C \to \M$ 
and certain associativity
and unit constraints satisfying some natural axioms (see \cite{O1} and 
references therein). 
The dual of a tensor category $\C$ with respect to an indecomposable
right module category is the category $\C^*_\M:=Fun_\C(\M,\M)$ 
whose objects are $\C$-module functors from $\M$ to itself and morphisms 
are natural module transformations. The category $\C^*_\M$ is a tensor category with 
tensor product being composition of module functors. 
Moreover, $\M$ becomes a left module category over $\C^*_\M$ in an obvious way.
The duality of tensor categories is known to be an equivalence relation \cite{Mu}.
A {\em fusion category} over $k$ is a 
$k$-linear semisimple rigid tensor category with finitely many 
isomorphism classes of simple objects and finite-dimensional Hom-spaces
such that the neutral object in simple (see \cite{ENO}). 

A finite tensor category
is said to be {\em pointed} if all its simple object are invertible. 
Every pointed semisimple category is equivalent to the fusion category 
$Vec(G, \, \omega)$ whose objects are vector spaces graded by the finite
group $G$ and whose associativity constraint is given by the $3$-cocycle 
$\omega \in Z^3(G, \, k^\times)$. Let us denote $Vec(G):= Vec(G, 1)$.
A fusion category is called {\em group-theoretical}
if it is equivalent to the dual of $Vec(G, \, \omega)$ with respect to
some indecomposable right module category for some finite group $G$ and
$3$-cocycle $\omega \in Z^3(G, \, k^\times)$.  

In this paper we use the notion of {\em weak Morita equivalence} 
\cite{Mu} of tensor categories
to define and study an equivalence relation called {\em categorical Morita equivalence} on the set
of all finite groups and on the set of all pairs $(G, \, \omega)$, where $G$ is a finite
group and $\omega \in H^3(G, \, k^\times)$. Namely, we say that two groups
$G$ and $G'$ (respectively, two pairs
$(G,\, \omega)$ and $(G',\, \omega')$)  are categorically Morita equivalent
if $Vec(G)$ is dual to $Vec(G')$ (respectively, $Vec(G, \, \omega)$ is dual
to $Vec(G', \, \omega')$) with respect to some indecomposable right module category.
This equivalence relation extends the notion of  {\em isocategorical groups},
i.e., groups with equivalent tensor categories of representations, studied in \cite{D}
and \cite{EG}. Our motivation to study categorical Morita equivalence of  finite groups
comes from the question about existence of semisimple Hopf
algebras with non group-theoretical representation categories 
asked in \cite[Question 8.45]{ENO}. We think that understanding equivalence 
classes of categorically  Morita equivalent groups is a natural step towards
answering this question.

The main results of this paper are: (1) Computation of the dual of 
$Vec(G, \, \omega)$ with respect to an indecomposable module category when
the dual is pointed, including explicit formulas for the Grothendieck
ring and the associated $3$-cocycle. 
(2) Necessary and sufficient conditions for two pointed semisimple categories 
to be dual to each other with respect to a module category.  
(3) A series of concrete examples of pairs of
groups $(G_1, \, G_2)$ that are categorically Morita equivalent but have 
non-isomorphic Grothendieck rings (and hence, inequivalent
representation categories). 
A consequence of the categorical Morita equivalence of these 
groups is that the representation categories of their Drinfeld doubles Rep(D($G_1$))
and Rep(D($G_2$)) are equivalent as braided tensor categories and so in particular 
these groups define
the same modular data. To the best of our knowledge  these are first examples of finite groups
with this property, cf.\ a discussion of a finite group modular data in \cite{CGR}.
These results are contained in
Theorems \ref{thm2}, \ref{thm3}, \ref{thm4} and Corollary \ref{thm5}.

The paper is organized as follows: in Section 2 we recall necessary
definition and facts from homological algebra. We also recall the notions
of module categories and duals of tensor categories. In Section 3
we give necessary and sufficient conditions for the dual of 
$Vec(G, \, \omega)$ with respect to an indecomposable right module category
to be pointed. In Section 4 we show that the
Grothendieck ring of the dual of $Vec(G, \, \omega)$ with respect to an 
indecomposable right module category when the dual is
pointed is the group ring of a certain crossed product
of groups. We also find an explicit formula for the 3-cocycle associated to the 
dual category. In Section 5 we introduce the
{\em categorical Morita equivalence} on the set
of all finite groups and on the set of all pairs $(G, \, \omega)$, where $G$ is a finite
group and $\omega \in H^3(G, \, k^\times)$. We give a 
group-theoretical and cohomological interpretation of these relations. 
In the final section, Section 6, we give a series of examples 
of pairs of groups that are categorically Morita equivalent 
but have non-isomorphic Grothendieck rings.

All categories considered in this paper are assumed to abelian, semisimple
and $k$-linear with finite dimensional Hom-spaces.
We will also assume that the number of isomorphism classes
of simple objects in any category is finite. All functors are assumed to be additive
and $k$-linear.

\end{section}


\begin{section}
{Preliminaries}

\begin{subsection}
{Cohomology of groups and Shapiro's lemma}

Let $G$ be a finite group and $M$ be a left $G$-module with action denoted by
$(g, \, m) \mapsto g \triangleright m$, for $g \in G, \, m \in M$. We define 
a cochain complex $C(G,\, M) = (C^n(G,\, M))_{n\geq 0}$
of $G$ with coefficients in $M$ as follows. Let $G^n = G \times \dots
\times G$ ($n$ factors) and $C^n(G,\, M) = Fun(G^n,\, M)$ be the set of all
$n$-cochains. By convention, $G^0(G, M) = M$. 
A $n$-cochain $f$ is said to be normalized if $f(g_1, \, g_2, \, \dots, g_n) = 0_M$
whenever $g_i = 1_G$ for some $i \in \{1, \, 2, \dots, n\}$.
All $n$-cochains are assumed to be normalized.
Let $\delta^n : C^n(G,\, M) \to C^{n+1}(G,\, M)$ be the coboundary operator
given by
\begin{eqnarray*}
(\delta^n f)(g_1, \dots, g_{n+1}) &=& g_1 \triangleright f(g_2,\dots, g_{n+1}) \\
& & +  \sum_{i=1}^n \, (-1)^i f(g_1,\dots, g_{i-1}, \, g_i g_{i+1}, \dots, g_{n+1}) \\
& & + (-1)^{n+1}f(g_1,\dots, g_{n})
\end{eqnarray*} 
for all $f\in C^n(G,\, M)$.  

If $M$ is a right $G$-module, we denote the action by
$(m, \, g) \mapsto m \triangleleft g$, for $g \in G, \, m \in M$.
Also, define $\underline{\delta}^n : C^n(G,\, M) \to C^{n+1}(G,\, M)$ by
\begin{eqnarray*}
(\underline{\delta}^n f)(g_1, \dots, g_{n+1}) &=& f(g_2,\dots, g_{n+1}) \\
& & +  \sum_{i=1}^n \, (-1)^i f(g_1,\dots, g_{i-1}, \, g_i g_{i+1}, \dots, g_{n+1}) \\
& & + (-1)^{n+1}(f(g_1,\dots, g_{n}) \triangleleft g_{n+1})
\end{eqnarray*} 
for all $f\in C^n(G,\, M)$.  

Let $Z^n(G,\, M) = Ker(\delta^n)$ be the set of $n$-cocycles and
$B^n(G,\, M) = Im(\delta^{n-1})$ be the space of $n$-coboundaries.
Let $\underline{Z}^n(G,\, M) = Ker(\underline{\delta}^n)$ and
$\underline{B}^n(G,\, M) = Im(\underline{\delta}^{n-1})$.
The $n$-th cohomology group $H^n(G,\, M)$ of $G$ with coefficients
in $M$ is the quotient $Z^n(G,\, M)/B^n(G,\, M)$ \\ $(n \geq 1)$. 
Also, let $\underline{H}^n(G,\, M) = \underline{Z}^n(G,\, M)/
\underline{B}^n(G,\, M)$.
\noindent When we write an element of the cohomology groups as 
$\overline{\omega}$, we will mean by this the class represented by the
cocycle $\omega$.  

Any homomorphism $a:G^\prime \to G$ between the groups $G$ and $G^\prime$
induces a map between their cohomology groups:

\begin{equation}
H^n(G,\, M) \to H^n(G^\prime,\, M) : 
(a, \overline{\omega})\mapsto  \omega^a := \omega \circ a^{\times n},
\end{equation}

Let $H$ be a subgroup of $G$.
Let $p:G \to H \setminus G$ be the usual surjection, i.e., $p(g) := Hg$,
for all $g \in G$. Throughout this paper we will denote $p(1_G)$ by $1$. 
For each $x \in H \setminus G$ choose a representative $u(x)$ in G; i.e., 
an element $u(x)$ with $pu(x) = x$. In particular, choose 
$u(1) = 1_G$. The set $H \setminus G$ is a right $G$-set
with the obvious action: $x \triangleleft g := p(u(x)g)$, for $x \in H \setminus G$ 
and $g \in G$. Also, the set $\{u(x) \, | \, x \in H \setminus G\}$ is a 
right $G$-set: $u(x) \triangleleft g = u(x \triangleleft g)$, for $x \in H \setminus G$ 
and $g \in G$ . The elements $u(x)g$ and
$u(x \triangleleft g)$ differ by an element $\kappa_{x, \, g}$ of $H$,
for $x \in H \setminus G$ and $g \in G$ :

\begin{equation}
\label{definition kappa}
u(x)g =  \kappa_{x,\, g}u(x \triangleleft g)
\end{equation}

\noindent The following relation holds:
\begin{equation}
\label{kappa relation}
\kappa_{x,\, g_1g_2} = \kappa_{x, \, g_1}\kappa_{x \triangleleft g_1, \, g_2}
\end{equation}
for any $x \in H \setminus G$ and $g_1, \, g_2 \in G$

The set $Fun(H \setminus G, \, k^\times)$ of functions from $H \setminus G$ to 
$k^\times$ is a left $G$-module: $(g \triangleright  f)(x) = f(x \triangleleft g)$,
for $x \in H \setminus G$ and $g \in G$ .
Let us regard $k^\times$ as a trivial left $H$-module. It is
easy to see that $Fun(H \setminus G, \, k^\times)$ is isomorphic to
the coinduced module $Coind_H^G(k^\times) = Hom_H(G, \,
k^\times)$. From now on we will identify the coinduced module
$Coind_H^G(k^\times)$ with $Fun(H \setminus G, \, k^\times)$.

Let $C := Coind_H^G(k^\times)$ and $K := H \setminus G$.
The action of $G$ on $K$ restricts to an action of $H$ on $K$.
Let $K^H$ denote the set of elements of $K$
that are stable under the action of $H$. Note that $K^H$ forms a group that
is isomorphic to $H \setminus N_G(H)$, where $N_G(H)$ is the normalizer 
of $H$ in $G$. By $\widehat{H}$, we will mean the group $Hom(H, \, k^\times)$
of 1-dimensional representations of $H$.

By Shapiro's Lemma there is an isomorphism between $H^n(G,\,
C)$ and $H^n(H,\, k^\times)$ for each
$n \in \mathbb{N}$. It is well known that the restriction maps
induces this isomorphism. We will need the explicit form of the
inverse of the restriction map when $n = 1, \, 2$. Lemmas 
\ref{shapiro1} and \ref{shapiro2} provide this.

\begin{lemma}
\label{shapiro1}
The following map induces an isomorphism between
$H^1(H,\, k^\times) = \widehat{H}$ and \\
$H^1(G,\, C)$:

\begin{equation}
\label{varphi_1}
\varphi_1:Z^1(H, \, k^\times) \to 
Z^1(G, \, C), \qquad (\varphi_1(\rho)(g))(x) = \rho (\kappa_{x,\, g})
\end{equation}

\noindent for any $\rho \in Z^1(H, \, k^\times)$, $g \in G$, $x \in K$.

\end{lemma}

\begin{proof} 
We will first show that $\varphi_1(\rho) \in
Z^1(G, \, C)$ for any $\rho \in
Z^1(H, \, C^\times)$. We need to show that $\varphi_1(\rho)$ satisfies the
equation:
\begin{equation*}
(\varphi_1(\rho)(g_1))(x) \, \,
(\varphi_1(\rho)(g_2))(x \triangleleft g_1)
= (\varphi_1(\rho)(g_1g_2))(x)
\end{equation*}
\begin{equation*}
\Leftrightarrow
\rho(\kappa_{x, \, g_1}) \rho(\kappa_{x \triangleleft g_1, \, g_2})
= \rho(\kappa_{x, \, g_1g_2})
\end{equation*}
The 1-cocycle condition on $\rho$ is:

\begin{equation*}
\rho(h_1) \rho(h_2) = \rho(h_1h_2)
\end{equation*}
Put $h_1 = \kappa_{x \, g_1}$ and $h_2 = \kappa_{x \triangleleft g_1, \, g_2}$ 
in the above equation and use Equation (\ref{kappa relation}) to
obtain the desired equation.

The map $\varphi_1$ induces a map:
\begin{equation*}
\widetilde{\varphi_1}:H^1(H, \, k^\times)
\to H^1(G, \, C)
\end{equation*}

\noindent One can show that 
the map $\widetilde{\varphi_1}$ does not depend
on the choice of the function $u: K \to G$.

\noindent Let $\psi_1$ denote the restriction map.
\begin{equation}
\label{psi_1}
\psi_1:Z^1(G, \, C) \to 
Z^1(H, \, k^\times), \qquad \psi_1(\gamma)(h) = \gamma(h)(1)
\end{equation}

\noindent for any $\gamma \in
Z^1(G, \, C)$ and $h\in H$.
Let $\widetilde{\psi_1}$ denote the induced map:

\begin{equation}
\widetilde{\psi_1}:H^1(G, \, C) \to
H^1(H, \, k^\times)
\end{equation}

It remains to show that the maps $\widetilde{\varphi_1}$ and
$\widetilde{\psi_1}$ are inverse to each other. It suffice to show that
$\widetilde{\psi_1} \circ \widetilde{\varphi_1}
= Id_{H^1(H, \, C^\times)}$. Actually, we show that
$\psi_1 \circ \varphi_1 = Id_{Z^1(H, \, k^\times)}$.
Pick any $\rho \in Z^1(H, \, k^\times)$. Then $\psi_1(\varphi_1(\rho))(h)
= (\varphi_1(\rho)(h))(1) 
= \rho(\kappa_{1, \, h}) = \rho(h)$ for all $h \in H$ and the lemma is
proved.

\end{proof}

\begin{lemma}
\label{shapiro2}
The following map induces an isomorphism between $H^2(H,\, k^\times)$ and \\
$H^2(G,\, C)$:

\begin{equation}
\label{varphi}
\varphi:Z^2(H, \, k^\times) \to 
Z^2(G, \, C), \qquad
(\varphi(\mu)(g_1,\, g_2))(x) = \mu (\kappa_{x,\, g_1}, \, 
\kappa_{x \triangleleft g_1, \, g_2}) \\
\end{equation}

for any $\mu \in Z^2(H, \, k^\times)$,
$g_1, \, g_2 \in G$, $x \in K$.

\end{lemma}

\begin{proof}
We will first show that $\varphi(\mu) \in
Z^2(G, \, C)$ for any $\mu \in
Z^2(H, \, k^\times)$. We need to show that $\varphi(\mu)$ satisfies the 
following equation for all $g_1, \, g_2, \, g_3 \in G$ and
$x \in K$.

\begin{equation*}
(\varphi(\mu)(g_2, \, g_3))(x \triangleleft g_1) \, \,
(\varphi(\mu)(g_1, \, g_2g_3))(x)
= (\varphi(\mu)(g_1g_2, \, g_3))(x) \, \,
(\varphi(\mu)(g_1, \, g_2))(x)
\end{equation*}
\begin{equation*}
\Leftrightarrow
\mu(\kappa_{x \triangleleft g_1, \, g_2}, \, \kappa_{x \triangleleft g_1g_2, \, g_3}) \, \, 
\mu(\kappa_{x, \, g_1}, \, \kappa_{x \triangleleft g_1, \, g_2g_3})
=\mu(\kappa_{x, \, g_1g_2}, \, \kappa_{x \triangleleft g_1g_2, \, g_3}) \, \,
\mu(\kappa_{x, \, g_1}, \, \kappa_{x \triangleleft g_1, \, g_2}) \\
\end{equation*}

\noindent The 2-cocycle condition on $\mu$ is:

\begin{equation*}
\mu(h_2, \, h_3) \, \mu(h_1, \, h_2h_3) = \mu(h_1h_2, \, h_3) \, \mu(h_1, \, h_2)
\end{equation*}

\noindent for any $h_1, \, h_2, \, h_3 \in H$.
Put $h_1 = \kappa_{x, \, g_1}, \, h_2 = \kappa_{x \triangleleft g_1, \, g_2}$ and 
$h_3 = \kappa_{x \triangleleft g_1g_2, \, g_3}$ in the above equation
and use equation (\ref{kappa relation}) to obtain the desired equation.

One can show that $\varphi$ preserves coboundaries, hence it induces a map:

\begin{equation}
\label{varphi tilde}
\widetilde{\varphi}:H^2(H, \, k^\times) \to
H^2(G, \, C)
\end{equation}

\noindent Let $\psi$ denote the restriction map:

\begin{equation}
\label{psi}
\psi:Z^2(G, \, C) 
\to Z^2(H, \, k^\times), \qquad
\psi(\gamma)(h_1, \, h_2) = \gamma(h_1, \, h_2)(1)
\end{equation}

\noindent for any $\gamma \in Z^2(G, \, C)$ 
and $h_1, \, h_2 \in H$. Let $\widetilde{\psi}$ denote the induced map:

\begin{equation}
\widetilde{\psi}: H^2(G, \, C) \to
H^2(H, \, k^\times) 
\end{equation}

It remains to show that the maps $\widetilde{\varphi}$ and $\widetilde{\psi}$
are inverse to each other. It suffice to show that
$\widetilde{\psi} \circ \widetilde{\varphi}
= Id_{H^2(H, \, k^\times)}$. Actually, we show that
$\psi \circ \varphi = Id_{Z^2(H, \, k^\times)}$.
Pick any $\mu \in Z^2(H, \, k^\times)$. Then $\psi(\varphi(\mu))
(h_1, \, h_2) = (\varphi(\mu)(h_1, \, h_2))(1) =
\mu(\kappa_{1, \, h_1}, \, \kappa_{1 \triangleleft h_1, \, h_2}) = 
\mu(h_1, \, h_2)$ for all $h_1, \, h_2 \in H$ and the lemma is proved.
\end{proof}

\noindent There is a right action of $K^H$ on $C^n(G, \, C)$:

\begin{equation*}
(\gamma, \, x) \mapsto {}^x \gamma, \qquad {}^x\gamma(g_1, \dots, g_n)(y) 
:= \gamma(g_1, \dots, g_n)(p(u(x)u(y)) \\
\end{equation*}

\vspace{0.2in}

\noindent for all $\gamma \in C^n(G, \, C), \, g_1, \dots, g_n \in G, \, 
x \in K^H$, and $y \in K$ \\

\noindent  It is routine to check that the above action is independent of 
the function $u$. This induces a right action of $K^H$ on $Z^n(G, \, C)$ 
and $H^n(G, \, C)$. If $H$ is normal in $G$, then $K^H = K$ and 

\begin{equation*}
{}^x\gamma(g_1, \dots, g_n)(y) := \gamma(g_1, \dots, g_n)(xy) \\ \\
\end{equation*}

\vspace{0.2in}

\noindent for all $\gamma \in C^n(G, \, C), \, g_1, \dots, g_n \in G, \, x, \, y \in K$. \\ 

\noindent Also, if $H$ is normal in $G$, then $Z^n(H, \, k^\times)$ is a right $G$-module:

\begin{equation*}
(\mu, \, g) \mapsto \mu^g, \qquad \mu^g(h_1, \dots, h_n) = \mu(gh_1g^{-1}, 
\dots, gh_ng^{-1}) \\ \\
\end{equation*}

\vspace{0.2in}

\noindent for all $\mu \in Z^n(H, \, k^\times), \, g \in G$ and 
$h_1, \dots, h_n \in H$. \\

\noindent If $H$ is abelian and normal in $G$, then $Z^n(H, \, k^\times)$ 
becomes a right $K$-module:

\begin{equation*}
(\mu, \, x) \mapsto \mu^{u(x)}
\end{equation*}

\vspace{0.2in}

\noindent for all $\mu \in Z^n(H, \, k^\times)$ and $x \in K$. \\

\noindent The above induces an action of $K$ on $H^n(H, \, k^\times)$.

\begin{lemma}
\label{hg-linear}
If $H$ is abelian and normal in $G$, then the map $\psi_1$ defined 
in (\ref{psi_1}) is a $K$-module map.
\end{lemma}

\begin{proof}
Pick any $\gamma \in Z^1(G, \, C)$ and $y \in K$.
We have $\psi_1({}^y\gamma)(h) = ({}^y\gamma)(h)(1)
= \gamma(h)(y)$ and
$(\psi_1(\gamma)^y)(h) = \psi_1(\gamma)(u(y)hu(y)^{-1})
= \gamma(u(y)hu(y)^{-1})(1)$. By Lemma \ref{shapiro1} we know that 
$\gamma = (\delta^1 \alpha) \, \varphi_1(\rho)$ for some $\alpha \in C$
and $\rho \in \widehat{H}$.

\noindent We have,

\begin{equation*}
\begin{split}
\gamma(h)(y) 
&= ((\delta^1 \alpha) \, \varphi_1(\rho))(h)(y) \\
&= \frac {\alpha(y \triangleleft h)} {\alpha(y)} \, \rho(\kappa_{y, \, h}) \\
&= \rho(u(y)hu(y)^{-1}) \\
\end{split}
\end{equation*}

\noindent and 

\begin{equation*}
\begin{split}
\gamma(u(y)hu(y)^{-1})(1) 
&= ((\delta^1 \alpha) \, \varphi_1(\rho))(u(y)hu(y)^{-1})(1) \\
&= \frac {\alpha(1 \triangleleft u(y)hu(y)^{-1})} {\alpha(1)} \, 
\rho(\kappa_{1, \, u(y)hu(y)^{-1}}) \\
&=\rho(u(y)hu(y)^{-1}). \\
\end{split}
\end{equation*}

\end{proof}

\begin{lemma}
\label{hg-linear1}
If $H$ is abelian and normal in $G$, then the map $\widetilde{\varphi}$ defined in
(\ref{varphi tilde}) is a $K$-module map. 
\end{lemma}

\begin{proof}
Pick any $\mu \in Z^2(H, \, k^\times)$. In order to show that the map 
$\widetilde{\varphi}$ is $K$-linear it suffices to show that 
$\psi({}^y\varphi(\mu))$ is cohomologous to $\psi(\varphi(\mu^y)) = \mu^y$ in $
H^2(H, \, k^\times)$ for all $y \in K$. We will actually show that 
$\psi({}^y\varphi(\mu))= \mu^y$. We have,

\begin{equation*}
\begin{split}
\psi({}^y\varphi(\mu))(h_1, \, h_2)
&= {}^y\varphi(\mu)(h_1, \, h_2)(1) \\
&= \varphi(\mu)(h_1, \, h_2)(y) \\
&= \mu(\kappa_{y, \, h_1}, \, \kappa_{y \triangleleft h_1, \, h_2}) \\
&= \mu(\kappa_{y, \, h_1}, \, \kappa_{y, \, h_2}) \\
&= \mu(u(y)h_1u(y)^{-1}, \, u(y)h_2u(y)^{-1}) \\
&= \mu^{y}(h_1, \, h_2) \\
\end{split}
\end{equation*}

\noindent for all $h_1, \, h_2 \in H$. So $\psi({}^y\varphi(\mu))= \mu^y$
and the lemma is proved.

\end{proof}

\end{subsection}

\begin{subsection}
{The fusion category {\bf $Vec(G, \, \omega)$}}

We refer the reader to \cite{BK} for definition and basic properties
of tensor categories.
A category is called {\em skeletal}
if all isomorphic objects in the category are actually equal. Every category 
is equivalent to a skeletal category. 
It is convenient to work with a skeletal category $\V(G,\, \omega)$
equivalent to $Vec(G, \, \omega)$. Let $\V(G,\,\omega)$ be
a semisimple tensor category with simple objects $g, \, g \in G$. The tensor product
is defined by $g_1 \otimes g_2 = g_1g_2$, and the associativity isomorphisms
are $\omega(g_1,\, g_2,\, g_3) id_{g_1g_2g_3}$. The unit object is $1_G$. 
The left and right unit isomorphisms are $\omega(1_G, 1_G, g) id_g$ and 
$\omega(g, 1_G, 1_G)id_g$, repectively. The previous statement follows
from the triangle axiom for tensor categories.
Since we can assume that all cocycles are normalized, the left and right unit 
isomorphisms are the identity morphisms. 
The left and right dual objects of $g$
are $g^* = {}^*g = g^{-1}$.
If $G^\prime$ is another group and $\omega^\prime \in Z^3(G^\prime, k^\times)$, 
then $\V(G, \, \omega) \cong \V(G^\prime, \, \omega^\prime)$
if and only if there is an isomorphism $a : G \to G^\prime$ such that
$\omega^\prime$ and $\omega^a$ are cohomologous. 

\end{subsection}

\begin{subsection}
{Module categories}

Recall some definitions from \cite{O1}:

\begin{definition} 
A right {\em module category} over a tensor category 
$(\C, \, \otimes, \, 1_{\C}, \, \alpha, \, \lambda, \, \rho)$ with
unit object $1_{\C}$, associativity constraint $\alpha$, left unit 
constraint $\lambda$, and right unit constraint $\rho$, is 
a category $\M$ together with an exact bifunctor $\otimes: \M \times \C \to 
\M$ and functorial associativity and unit isomorphisms:
$\mu_{M, \, X, \, Y}: M \otimes (X \otimes Y) \to (M \otimes X) \otimes Y, \,\,
\tau_M: M \otimes 1_{\C} \to M$ for any $X, Y \in \C, \, M \in \M$ such that 
the following two diagrams commute.

\begin{equation}
\label{module pentagon}
\xymatrix{&M\otimes (X \otimes Y ) \otimes Z)
\ar[dl]_{\id_M \otimes \alpha_{X,Y,Z}} 
\ar[dr]^{\mu_{M, \, X \otimes Y, \, Z}}&\\
M \otimes (X \otimes (Y \otimes Z))
\ar[d]^{\mu_{M, \, X, \,  Y \otimes Z}}&&
(M \otimes (X \otimes Y)) \otimes Z 
\ar[d]_{\mu_{M, \, X, \, Y}\otimes \id_Z}\\
(M \otimes X) \otimes (Y \otimes Z) 
\ar[rr]^{\mu_{M \otimes X, \, Y, \, Z}}&&
((M \otimes X)\otimes Y) \otimes Z}
\end{equation}

\begin{equation}
\label{module triangle}
\xymatrix{M \otimes (1_\C \otimes Y)\ar[rr]^{\mu_{M, \, 1_\C, \, Y}} 
\ar[dr]^{\id_M \otimes \lambda_Y}&&(M \otimes 1_\C) \otimes Y
\ar[dl]_{\tau_M \otimes \id_Y}\\ &M \otimes Y&}
\end{equation}

\end{definition}

\begin{definition}
Let $(\M_1, \, \mu^1, \tau^1)$ and $(\M_2, \, \mu^2, \tau^2)$ be two
right module categories over
a tensor category $\C$. A {\em module functor} from $\M_1$ to $\M_2$
is a functor $F: \M_1\to \M_2$ together with functorial isomorphisms
$\gamma_{M, \, X}: F(M \otimes X) \to F(M) \otimes X$ for any
$X \in \C, \, M \in \M_1$ such that the following two diagrams commute.

\begin{equation}
\label{module functor pentagon}
\xymatrix{&F(M \otimes (X\otimes Y))\ar[dl]_{F(\mu^1_{M, \, X, \, Y})} \ar[dr]^{\gamma_{M, X\otimes 
Y}}&\\ F((M\otimes X)\otimes Y)\ar[d]^{\gamma_{M\otimes X, \, Y}}&&F(M)\otimes (X\otimes Y)
\ar[d]_{\mu^2_{F(M), \, X, \, Y}}\\ F(M\otimes X) \otimes Y \ar[rr]^{\gamma_{M, \, X} \otimes \id_Y}
&& (F(M)\otimes X)\otimes Y}
\end{equation}

\begin{equation}
\label{module functor triangle}
\xymatrix{F(M \otimes 1_\C)\ar[rr]^{F(\tau^1_M)} \ar[dr]^{\gamma_{M, \, 1_\C}}&&F(M)\\ &F(M) 
\otimes 1_\C\ar[ur]^{\tau^1_{F(M)}}&}
\end{equation}

\end{definition}

Two module categories $\M_1$ and $\M_2$ over $\C$ are {\em equivalent}
if there exists a module functor from $\M_1$ to $\M_2$ which is an
equivalence of categories.
For two module categories $\M_1$ and $\M_2$ over a tensor category
$\C$ their {\em direct sum} is the category $\M_1 \oplus \M_2$ with 
the obvious module category structure. A module category is 
{\em indecomposable} if it is not equivalent to a
direct sum of two non-trivial module categories.

\begin{definition}
Let $\M_1$ and $\M_2$ be two right module categories over
a tensor category $\C$. Let $(F^1, \, \gamma^1)$ and $(F^2, \, \gamma^2)$
be module functors from $\M_1$ to $\M_2$. A {\em natural module
transformation} from $(F^1, \, \gamma^1)$ to $(F^2, \, \gamma^2)$ is a 
natural transformation $\eta: F^1 \to F^2$ such that the following
square commutes for all $M \in \M$, $X \in \C$.

\begin{equation}
\label{module trans square}
\xymatrix{F^1(M \otimes X)\ar[r]^{\eta_{M \otimes X}}\ar[d]_{\gamma_{M, \, X}^1} 
&F^2(M \otimes X)\ar[d]^{\gamma_{M, \, X}^2}&\\
F^1(M) \otimes X \ar[r]_{\eta_M \otimes id_X} &F^2(M) \otimes X&} 
\end{equation}

\end{definition}

\begin{example}
Let us recall a description of indecomposable module categories 
over $\V(G, \, \omega)$ given in \cite{O2}. 
Let $\M$ be an indecomposable right module category
over $\V(G, \, \omega)$ with module category structure $\mu$.
Without loss of generality we may assume that
$\M$ is skeletal. The set of simple objects of $\M$ is a transitive right
$G$-set and hence can be identified with the set of right cosets
$H \setminus G = K$ for some subgroup $H$ of $G$.
So the set of all simple objects of $\M$,
$Irr(\M)$ = $K$. All the isomorphisms $\mu_{x, \, g_1, \, g_2}$,
$x \in K$, $g_1, g_2 \in G$ are given by scalars. So we can 
regard $\mu$ as an element of $C^2(G, \, C)$:

\begin{equation*}
\mu(g_1, \, g_2)(x) := \mu_{x, \, g_1, \, g_2}, 
\qquad x \in K, g_1, g_2 \in G
\end{equation*}

\noindent We may assume that the 
$2$-cochain $\mu$ is normalized. Since the unit constraint in 
$\V(G, \, \omega)$ is trivial, the commutativiy of triangle 
(\ref{module triangle}) implies that the unit constraint in $\M$ is trivial.
Let us regard $\omega$ as an element of $Z^3(G, \, C) \subset C^3(G, \, C)$ by
treating $\omega(g_1, \, g_2, \, g_3)$ as a constant function on $K$, for all
$g_1, \, g_2, \, g_3 \in G$. 
The commutativity of the pentagon (\ref{module pentagon}) implies that

\begin{equation}
\label{eq1}
\delta^2\mu = \omega 
\end{equation}

This in particular means that $\omega$ restricted
to $H \times H \times H$ represents the trivial class in $H^3(H, \, k^\times)$. 
Let $L_{H, \, \omega} := \{\mu \in C^2(G, \, C) \, | \,
\delta^2\mu = \omega \}$.
Two elements in $L_{H, \, \omega}$ give rise to equivalent module categories
if and only if the differ by some element in $B^2(G, \, C)$.
Define an equivalence relation on $L_{H, \, \omega}$: two elements in 
$L_{H, \, \omega}$ are equivalent if and only if the differ by an element
in $B^2(G, \, C)$.
We denote the set of equivalence classes of $L_{H, \, \omega}$ under the previous
relation by $\overline{L}_{H, \, \omega}$.
The sets $\overline{L}_{H, \, \omega}$ and $H^2(H, \, k^\times)$
are in bijection.

\end{example}

\end{subsection}


\begin{subsection}
{The dual category}

Let $\C$ be a tensor category and $\M$ be an indecomposable 
right module category over $\C$. 

\begin{definition} 
\label{dual category}
The dual category of $\C$ with respect to
$\M$ is  the category $\C^*_\M:=Fun_\C(\M,\M)$ whose objects are $\C$-module 
functors from $\M$ to itself and morphisms are natural 
module transformations.
\end{definition}

The category $\C^*_\M$ is tensor with tensor product being composition
of module functors.
Let $(\gamma^1, \, F^1)$, $(\gamma^2, \, F^2) \in Obj(\C^*_\M)$, where
$\gamma^1, \, \gamma^2$ represent the module functor structure on the
functors $F^1$ and $F^2$, respectively. Then,
$(\gamma^1, \, F^1) \otimes (\gamma^2, \, F^2) = (\gamma, \, F^1 \circ F^2)$,
where $\gamma$ is defined as: $\gamma_{M, \, X} := 
\gamma^1_{F^2(M), \, X} \circ F^1(\gamma^2_{M, \, X})$ for any 
$M \in \M$, $X \in \C$. 
Let $\eta: (\gamma^1, \, F^1) \to (\gamma^2, \, F^2)$ and 
$\eta^\prime: (\gamma^3, \, F^3) \to (\gamma^4, \, F^4)$ be morphisms in
$\C^*_\M$, i.e., natural module transformations. Then their tensor product
$\eta \otimes \eta^\prime$ is defined as: 
$(\eta \otimes \eta^\prime)(M) := \eta_{F^4(M)} \circ F^1(\eta^\prime_M)$.

\begin{remark} 
\label{dimension of dual}
The {\em Frobenius-Perron}
dimension $FPdim(X)$ of a simple object $X \in Obj(\C)$ is the 
{\em Frobenius-Perron eigenvalue} of the matrix coming from 
multiplication of the set of isomorphism classes of all simple objects in $\C$
by $X$. The Frobenius-Perron dimension $FPdim(\C)$ of the fusion category $\C$
is the sum of squares of the Frobenius-Perron dimension of the isomorphism classes
of simple objects.
It is know that if $\C$ is a fusion category and $\M$ is abelian and semisimple 
then $\C^*_\M$ is a fusion category. 
It is also known that $FPdim(\C) = FPdim(\C^*_\M)$. See \cite{ENO} for a treatment
on Frobenius-Perron dimensions.
\end{remark}

\end{subsection}


\end{section}

\begin{section}
{Necessary and sufficient condition for the dual of a pointed category to be pointed}

We fix the following notation for this and the next section.
Let $K := H \setminus G$ and $C:= Coind^G_H k^\times$. Let $u: K \to G$ 
be a function satisfying 
$p \circ u = id_{K}$ and $u(p(1_G)) = 1_G$, where $p:G \to K$
is the usual surjection. Let $\kappa:
K \times G \to H$ be the function satisfying Equation (\ref{definition kappa}).
Let $\C := \V(G, \, \omega)$ and let $\M = \M(H, \, \mu)$ 
denote the right module category constructed from the pair $(H,\, \mu)$, 
where $H$ is a subgroup of $G$ such that $\omega|_{H \times H \times H}$ is
trivial in $H^3(H, \, k^\times)$
and $\mu \in C^2(G, \, C)$ is a 2-cochain satisfying $\delta^2\mu = \omega$.
In the previous equation we regarded $\omega$ as an element of $Z^3(G, \, C)$ 
by treating $\omega(g_1, \, g_2, \, g_3)$ as a constant function on $K$, for all
$g_1, \, g_2, \, g_3 \in G$. 
The module category structure of $\M$ is given by $\mu$.
If $\omega \equiv 1$, then we will assume that $\mu$ belongs to 
$Z^2(H, \, k^\times)$ and that the module category structure of 
$\M(H, \, \mu)$ is given by $\varphi(\mu)$ (see (\ref{varphi})).

\begin{definition}
For each $y \in K^H$, define the set
$Fun_y = Fun_y(G, \, C)$:

\begin{equation*}
Fun_y := \left\{\gamma \in C^1(G, \, C) \, | \, \delta^1 \gamma
= \frac {{}^y\mu} {\mu} \right\}
\end{equation*}
\end{definition}

\begin{lemma}
\label{invertible objects}
Invertible objects in $\C_\M^*$ are given by pairs $(\gamma, \, y)$, where
$y \in K^H$ and $\gamma \in Fun_y$.
\end{lemma}

\begin{proof}
We associate an invertible objects in $\C_\M^*$ to each pair $(\gamma, \, y)$, 
where $y \in K^H$ and $\gamma \in Fun_y$ as follows: define a map 
$f_y: K \to K$ by $f_y(x) = p(u(y)u(x))$ for any $x \in K$. Extend the map $f_y$ 
to a functor  $F_y: \M \to \M$. 
The module functor structure on $F_y$, which is also denoted by $\gamma$, 
is: $\gamma_{x, \, g} := \gamma(g)(x) \, id_{p(u(y)u(x \triangleleft g))}$ 
for any $g \in G $ and $x \in K$. The pentagon axiom 
for a module functor (\ref{module functor pentagon}) is:

\begin{equation*}
{}^y\mu(g_1, \, g_2)(x) \,\, \gamma(g_1g_2)(x) =
\gamma(g_1)(x) \, \, \gamma(g_2)(x \triangleleft g_1) \,\, 
\mu(g_1, \, g_2)(x)
\end{equation*} 
for all $g_1, \, g_2 \in G$ and $x \in K$.

This condition is satisfied because $\gamma \in Fun_y$.
The inverse of $(\gamma, \, F_y)$ is the module functor associated to the pair
$({}^{p(u(y)^{-1})}\gamma^{-1}, \, p(u(y)^{-1}))$,
where ${}^{p(u(y)^{-1})}\gamma^{-1}$ is defined as: \\
${}^{p(u(y)^{-1})}\gamma^{-1}(g)(x) := {\gamma(g)(p(u(y)^{-1}u(x))}^{-1}$, 
for $g \in G$ and $y \in K$. It should be clear that all
invertible objects in $\C_\M^*$ arise in this way and the lemma is proved.
\end{proof}

Two invertible $\C$-module functors $(\gamma^1, \, y_1)$  and
$(\gamma^2, \, y_2)$ are isomorphic in ${\C}_{\M}^{*}$ iff $y_1 = y_2$
and there exists an element $\alpha \in C$
such that $\gamma^1(g)(x) = \frac {\alpha(x \triangleleft g)} {\alpha(x)}
\gamma^2(g)(x)$ for all $g \in G$ and $x \in K$.

This motivates us to define an equivalence relation on the set
$Fun_y$: we define two elements $\gamma^1, \, \gamma^2 \in Fun_y$,
to be equivalent if there exists an $\alpha \in C$
such that 

\begin{equation*}
\gamma^1(g)(x) = \frac {\alpha(x \triangleleft g)} {\alpha(x)}
\gamma^2(g)(x)
\end{equation*}

\noindent for all $g \in G$ and $x \in K$.

\noindent Let $\overline{Fun}_y$ denote the set of equivalence classs
of $Fun_y$ under the aforementioned equivalence relation.

\begin{lemma}
\label{al}
For each $y \in K^H$, if $\overline{Fun}_y \neq \emptyset$, then 
there is a bijection between the sets $\overline{Fun}_y$ and 
$H^1(G, \, C)$ and hence there is a bijection between the
sets $\overline{Fun}_y$ and $\widehat{H}$.
\end{lemma}

\begin{proof}
Fix some $\eta_y \in Fun_y$. Then the maps $Fun_y \to 
Z^1(G, \, C): \beta \mapsto \frac {\beta} {\eta_y}$
and $Z^1(G, \, C) \to Fun_y: \gamma \mapsto \eta_y \gamma$
are inverse to each other. These maps induce a bijection
between the sets $\overline{Fun}_y$ and $H^1(G, \, C)$
The second statement of the lemma follows from
Shapiro's Lemma.
\end{proof}

\begin{theorem}
\label{thm1}
The tensor category ${\C}_{\M}^{*}$ (where $\M$ is the 
$\C = V(G, \, \omega)$-module category constructed from the 
pair $(H, \, \mu)$ where 
$H$ is a subgroup of $G$ such that $\omega|_{H \times H \times H}$ is
trivial in $H^3(H, \, k^\times)$ and $\mu \in C^2(G, \, C)$ is a 2-cochain 
satisfying $\delta^2\mu = \omega$) is pointed if and only if
the following three conditions hold:
\begin{enumerate}
\item $H$ is abelian,
\item $H$ is normal in $G$ and
\item the restriction $\psi({}^y\mu/\mu)$ is trivial in $H^2(H, \, k^\times)$, for
all $y \in K$.
\end{enumerate}

If $\omega \equiv 1$, then we assume that $\mu$ belongs to $Z^2(H, \, k^\times)$
and the module category structure on $\M$ is given by $\varphi(\mu)$ (see (\ref{varphi})).
The third condition above is then replaced with:

\begin{enumerate}
\item[(3${}^\prime$)]

$\mu$ represents a $G$-invariant class in $H^2(H, \, k^\times)$.
\end{enumerate}

\end{theorem}

\begin{proof}
Suppose that ${\C}_{\M}^{*}$ is pointed and
let $\S = K^H$. The set of isomorphism classes of simple objects in 
${\C}_{\M}^{*}$ is given by the set $\bigcup_{s \in \S}
\left(\overline{Fun_s} \times \{s\}  \right)$.
By the previous Lemma, we have $FPdim({\C}_{\M}^{*}) \leq |\widehat{H}| \, |\S|$.
Note that $|\widehat{H}|\leq |H|$ and $|\S| \leq |K| =
\frac {|G|} {|H|}$.
By Remark \ref{dimension of dual}, $FPdim({\C}_{\M}^{*}) = FPdim(\C) = |G|$. 
It follows that we must have $\overline{Fun_y} \neq \emptyset$ for all $y \in K$,
$|\widehat{H}| = |H|$ and $\S = K$. The second condition in the previous sentence
means that $H$ is abelian. The third condition means that $H$ is normal in $G$. 
The first condition is equivalent to saying that
$\frac{{}^y\mu}{\mu}$ is trivial in $H^2(G, \, C)$, for
all $y \in K$. This is equivalent to saying that the restriction 
$\psi(\frac{{}^y\mu}{\mu})$ is trivial in $H^2(H, \, k^\times)$,
for all $y \in K$.

Conversely, suppose that $H$ is abelian and normal in $G$ and
that $\psi(\frac{{}^y\mu}{\mu})$ is trivial in $H^2(H, \, k^\times)$, for
all $y \in K$.
Let $\C^\prime$ denote the full subcategory of $\C_\M^*$
of invertible objects. The isomorphism classes of invertible objects in 
the category ${\C}_{\M}^{*}$ are given by elements of the set
$\bigcup_{x \in K} \left(\overline{Fun_x} \times \{x\}\right)$ and the
size of each set in the previous union is $|H|$. 
So $FPdim(\C^\prime) = |G|$. It follows that  ${\C}_{\M}^{*} = \C^\prime$.
In other words, every simple object in ${\C}_{\M}^{*}$ is invertible, that is,
the category ${\C}_{\M}^{*}$ is pointed.

The last statement of the theorem follows from Lemma \ref{hg-linear}.
\end{proof}

\begin{example}
If $G = \mathbb{Z}/n\mathbb{Z}$ is a finite cyclic group, 
then $H^2(H, \, k^\times) = \{1\}$ for any subgroup $H$ of $G$. 
Hence the dual of $Vec(G, \, \omega)$
with respect to any indecomposable module category for any 3-cocycle 
$\omega$ on $G$ is pointed. Also,
for any abelian group $G$, the dual of $Vec(G)$ with respect to any
indecomposable module category is pointed. On the other hand, the previous
statement is not true for $Vec(G, \, \omega)$ if $\omega$ is a non-trivial 
3-cocycle on the abelian group $G$. Indeed, consider the dihedral group $D_8 = 
\{r, \, s \, | \, r^4 = s^2 = 1, \, rs = sr^{-1}\}$ and a subgroup $<~r>$ of it.
It can be shown that $Vec(D_8)^*_{\M(<r>, \, 1)} \cong 
Vec((\mathbb{Z}/2\mathbb{Z})^3, \, \omega)$, where $\omega$ is a non-trivial
3-cocycle on $(\mathbb{Z}/2\mathbb{Z})^3$. Now, we know that $Vec(D_8)$ is 
dual to the representation category
Rep($D_8$). Hence, there must exist an indecomposable module category over
$Vec((\mathbb{Z}/2\mathbb{Z})^3, \, \omega)$ with respect to which the dual
of $Vec((\mathbb{Z}/2\mathbb{Z})^3, \, \omega)$ is equivalent to the non-pointed
tensor category Rep($D_8$). We refer the reader to \cite{CGR} and \cite{GMN}
for similar results.
\end{example}

\end{section}


\begin{section}
{The tensor category ${\C}_{\M}^{*}$ when it is pointed}

In this section we further assume that $H$ is abelian and normal in $G$ and that
$\frac{{}^y\mu}{\mu}$ is trivial in $H^2(G, \, C)$, for all $y \in K$.


\begin{subsection}
{Tensor product and composition of morphisms in ${\C}_{\M}^{*}$}

It suffices to restrict ourselves to simple objects in ${\C}_{\M}^{*}$.
Recall that simple objects in ${\C}_{\M}^{*}$ are given by 
pairs $(\gamma, \, x)$, where $\gamma \in Fun_x$ and $x \in 
K$. The element $x \in K$ determines a $\C$-module functor
$F_x: \M \to \M$ given by $F_x(y) = xy$, for any 
$y \in K$. The $\C$-module functor structure on $F_x$ is
given by $\gamma$. Tensor product (=composition of module functors) in
${\C}_{\M}^{*}$: for any two simple objects $(\gamma^1, \, x_1)$ and $(\gamma^2, \, x_2)$, 
$(\gamma^1, \, x_1) \otimes (\gamma^2, \, x_2) =
({}^{x_2}\gamma^1 \, \gamma^2, \, x_1x_2)$ where ${}^{x_2}\gamma^1 \,
\gamma^2$ is an element of the set $Fun_{x_1x_2}$ and ${}^{x_2}\gamma^1$ is 
defined as follows: ${}^{x_2}\gamma^1(g)(y) = \gamma^1(g)(x_2y)$, for $g \in G$,
$y \in K$.

Now let us look at morphisms in ${\C}_{\M}^{*}$.
It suffices to restrict ourselves to isomorphisms between
simple objects. Recall that an isomorphism between two simple objects 
$(\gamma^1, \, x)$ and $(\gamma^2, \, x)$ (note that the second coordinates
have to be equal for an isomorphism to exist) in ${\C}_{\M}^{*}$ is given
by an element $\alpha \in C$ which satisfies:
$\gamma^1(g)(y) = \frac {\alpha(y \triangleleft g)} {\alpha(y)} \,
\gamma^2(g)(y)$, for all $g \in G$ and $y \in K$.

\begin{note}
\label{auto}
An isomorphism $\alpha : (\gamma^1, \, x) \to (\gamma^2, \, x)$  
is completely determined by $\alpha(1)$. If $\alpha$ is an automorphism, then
$\alpha(y) = \alpha(1)$ for all $y \in K$.
\end{note}

Now let us look at  tensor product and composition of isomorphisms
in ${\C}_{\M}^{*}$. Let 
$\alpha : (\gamma^1, \, x_1) \to (\gamma^2, \, x_1)$ and 
$\beta : (\gamma^3, \, x_2)\to (\gamma^4, \, x_2)$ be any two isomorphisms between 
simple objects in ${\C}_{\M}^{*}$. The tensor product of $\alpha$ and $\beta$: 
$(\alpha \otimes \beta)(x) = ({}^{x_2}\alpha \, \beta)(x)
= \alpha(x_2x) \beta(x)$ for any
$x \in K$. If $\gamma^2 = \gamma^3$, then the composition of 
$\alpha$ and $\beta$ is given by $(\beta \circ \alpha)(x) = \beta(x) \alpha(x)$ 
for $x \in K$. 

\end{subsection}


\begin{subsection}
{Grothendieck ring of the category ${\C}_{\M}^{*}$}

The set of isomorphism classes of simple objects in ${\C}_{\M}^{*}$ form a group:

\begin{equation}
\Lambda = \bigcup_{x \in K} \left(\overline{Fun}_x \times \{x\} \right) \qquad
(\overline{\gamma^1}, \, x_1) \star (\overline{\gamma^2}, \, x_2)
= (\overline{{}^{x_2}\gamma^1 \, \gamma^2}, \, x_1x_2)
\end{equation}

\noindent where for any $\gamma \in Fun_x$, by $\overline{\gamma}$ we mean the
equivalence class of $\gamma$ in $\overline{Fun}_x$. The inverse of any 
$(\overline{\gamma}, \, x) \in \Lambda$ is 
$\left(\overline{\gamma_{x^{-1}}^{-1}}, \, x^{-1}\right)$. 
The Grothendieck ring $\mathcal{K}_0({\C}_{\M}^{*})$ equals $\mathbb {Z}[\Lambda]$.

The rest of the subsection is devoted to showing that $\Lambda$ is 
isomorphic to a certain crossed product of the groups $\widehat{H}$ and
$K$.

Since we assumed that $\frac{{}^y\mu}{\mu}$ is trivial in $H^2(G, \, C)$,
for each $y \in K$ we have a map $\eta_y \in C^1(G, \, C)$ such that:

\begin{equation}
\label{eta1}
\delta^1 \eta_y = \frac {{}^y\mu} {\mu}
\end{equation}

\noindent Define a function 

\begin{equation}
\label{nutilde}
\tilde{\nu}: K \times K \to
C^1(G, \, C), \qquad
\tilde{\nu}(y_1, \, y_2) = \frac {{}^{y_2}\eta_{y_1} \, \eta_{y_2}}
{\eta_{y_1y_2}}
\end{equation}

\begin{lemma}
\label{nutilde1}
The function $\tilde{\nu}$ defines an element in
$\underline{H}^2(K, \, H^1(G, \, C))$.
\end{lemma}

\begin{proof}
Let us first show that $\tilde{\nu}(y_1, \, y_2) \in
Z^1(G, \, C)$ for any $y_1, \, y_2 
\in K$.
We have $\delta \eta_{y_1y_2}
= \frac {{}^{y_1y_2}\mu} {\mu}
= \frac {{}^{y_2}({}^{y_1}\mu)} {\mu}
= \frac {{}^{y_2}(\delta\eta_{y_1}\mu)} {\mu}
= \frac {{}^{y_2}(\delta\eta_{y_1}) \, {}^{y_2}\mu} {\mu}
= \delta({}^{y_2}\eta_{y_1}) \,\, \delta\eta_{y_2}
= \delta({}^{y_2}\eta_{y_1} \, \eta_{y_2})$.
So $\tilde{\nu}(y_1, \, y_2) \in Z^1(G, \, C)$ for any $y_1, \, y_2 \in K$.
Now let us show that $\tilde{\nu} \in Z^2(K, \, Z^1(G, \, C))$.
We have
\begin{equation*}
\begin{split}
(\delta\tilde{\nu})(y_1, \, y_2, \, y_3)
&=\tilde{\nu}(y_2, \, y_3) \tilde{\nu}(y_1y_2, \, y_3)^{-1}
\tilde{\nu}(y_1, \, y_2y_3) ({}^{y_3}\tilde{\nu}(y_1, \, y_2))^{-1} \\
&= \frac {{}^{y_3}\eta_{y_2} \, \eta_{y_3}}
{\eta_{y_2y_3}} \,\,\times\,\,
\frac {\eta_{y_1y_2y_3}}
{{}^{y_3}\eta_{y_1y_2} \, \eta_{y_3}} \,\,\times\,\,
\frac {{}^{y_2y_3}\eta_{y_1} \, \eta_{y_2y_3}}
{\eta_{y_1y_2y_3}} \,\,\times\,\,
\frac {{}^{y_3}\eta_{y_1y_2}}
{{}^{y_3}({}^{y_2}\eta_{y_1}) \, {}^{y_3}\eta_{y_2}} \\
&\equiv 1 \\
\end{split}
\end{equation*}

The cohomology class of $\tilde{\nu}$ does not depend on the choice of the 
family of maps $\{\eta_y \, | \, y \in K\}$. 
Indeed, let $\{\eta^\prime_y \, | \, y \in K\}$ be another family of maps
satisfying $(\delta \eta^\prime_y) = \frac {{}^y\mu} {\mu}$ for all 
$y \in K$. We want to show that $\tilde{\nu}(y_1, \, y_2) = 
\frac {{}^{y_2}\eta_{y_1} \, \eta_{y_2}} {\eta_{y_1y_2}}$ 
and $\tilde{\nu}^\prime(y_1, \, y_2) = \frac {{}^{y_2}\eta^\prime_{y_1} \, \eta^\prime_{y_2}}
{\eta^\prime_{y_1y_2}}$ define the same class in 
$\underline{H}^2(K, \, Z^1(G, \, C))$.
We have, $\delta(\frac {\eta_y} {\eta^\prime_y}) = 1$, i.e. 
$\frac {\eta_y} {\eta^\prime_y} \in Z^1(G, \, C)$
for each $y \in K$.
Define $\beta: K \to Z^1(G, \, C)$ by
$\beta(y) := \frac {\eta_y} {\eta^\prime_y}$. Then, 
$\tilde{\nu}(y_1, \, y_2) = 
\frac {{}^{y_2}\eta_{y_1} \, \eta_{y_2}} {\eta_{y_1y_2}}
= \frac {{}^{y_2}\beta(y_1)\, {}^{y_2}\eta^\prime_{y_1} \,\, \beta(y_2) \, \eta^\prime_{y_2}}
{\beta(y_1y_2) \eta^\prime_{y_1y_2}}
= (\delta \beta)(y_1, \, y_2) \,\, \tilde{\nu}^\prime(y_1, \, y_2)$.
\end{proof}

\begin{corollary}
\label{asso}
The function $\nu = \psi_1 \circ \tilde{\nu}$ defines an element in
$\underline{H}^2(K, \, \widehat{H})$.
\end{corollary}

\begin{proof}
The proof follows immediately from  Lemmas \ref{hg-linear} and \ref{nutilde1}.
\end{proof}

\begin{remark}

If $\omega \equiv 1$, then the element $\nu$ in the previous corollary is 
the image of $\mu$ under the following composition.

\begin{equation}
\label{Phi}
\Phi:H^2(H, \, k^\times)^{K}
\longrightarrow H^2(G, \, C)^{K}
\longrightarrow \underline{H}^2(K, \, H^1(G, \, C))
\longrightarrow \underline{H}^2(K, \, \widehat{H})
\end{equation}

The first map in the above composition comes from $\varphi$ (\ref{varphi}), 
the second from (\ref{nutilde}) and third is induced from the
map $\psi_1$ (\ref{psi_1}). Maps similar to 
the one in (\ref{Phi}) appears in \cite{D} and \cite{EG}.

\end{remark}

Let us put a group structure on the set $\widehat{H} \times K$.
For any two pairs $(\rho_1, \, x_1), \, (\rho_2, \, x_2)$ define their
product by:

\begin{equation}
\label{cproduct1}
(\rho_1, \, x_1)  (\rho_2, \, x_2) =
(\nu(x_1, \, x_2) \, \rho_1^{x_2} \rho_2, \, x_1x_2)
\end{equation}

\noindent Associativity follows from corollary \ref{asso}. We denote this group by
$\widehat{H} \rtimes_{\nu} K$. The group that we just constructed is known as a crossed
product.

As mentioned in Lemma \ref{al}, the sets $\overline{Fun}_x$ and $\widehat{H}$ are in
bijection for each $x \in K$. The following maps induce this bijection:

\begin{equation}
\label{zeta}
\zeta_x: \widehat{H} \to Fun_x, \qquad \zeta_x(\rho) := \eta_x \, \varphi_1(\rho) 
\end{equation}

\begin{equation*}
\label{theta}
\theta_x: Fun_x \to \widehat{H}, \qquad \theta_x(\gamma) := 
\psi_1(\gamma / \eta_x)
\end{equation*}

\noindent where the maps $\varphi_1$ and $\psi_1$ were defined in (\ref{varphi_1})
and (\ref{psi_1}), respectively.

\begin{theorem}
\label{thm2}
The Grothendieck ring $\mathcal{K}_0({\C}_{\M}^{*}) = \mathbb{Z}[\Lambda]$
is isomorphic to the group ring $\mathbb{Z} [\widehat{H} \rtimes_{\nu} 
K]$.
\end{theorem}

\begin{proof}
Suffices to show that the groups $\Lambda$ and $\widehat{H} \rtimes_{\nu} 
K$ are isomorphic. Define a map $T:\widehat{H} \rtimes_{\nu} 
K \to \Lambda$ by $T((\rho, \, x)) = (\overline{\zeta_x(\rho)}, \, x)$. 
Let us show that $T$ is a group homomorphism. For any $(\rho_1, \, x_1)$,
$(\rho_2, \, x_2) \in \widehat{H} \rtimes_{\nu} K$, we have 

\begin{equation*}
\begin{split}
\qquad \quad \quad \, T((\rho_1, \, x_1)(\rho_2, \, x_2))
&= T((\nu(x_1, \, x_2) \, \rho_1^{x_2} \rho_2, \, x_1x_2)) \\
&= (\overline{\zeta_{x_1x_2}(\nu(x_1, \, x_2) \, \rho_1^{x_2} \rho_2)}, \, x_1x_2) \\
\end{split}
\end{equation*}

and,

\begin{equation*}
\begin{split}
T((\rho_1, \, x_1)) \star T((\rho_2, \, x_2))
&= (\overline{\zeta_{x_1}(\rho_1)}, \, x_1) \star (\overline{\zeta_{x_2}(\rho_2)}, \, x_2) \\
&= (\overline{{}^{x_2}(\zeta_{x_1}(\rho_1)) \, \zeta_{x_2}(\rho_2)}, \, x_1x_2) \\
\end{split}
\end{equation*}

\noindent Now, we show that $\theta_{x_1x_2}({}^{x_2}\zeta_{x_1}(\rho_1) \, \zeta_{x_2}(\rho_2))
= \nu(x_1, \, x_2) \, \rho_1^{x_2} \rho_2$. For any $h \in H$, we have

\begin{equation*}
\begin{split}
\theta_{x_1x_2}({}^{x_2}(\zeta_{x_1}(\rho_1)) \, \zeta_{x_2}(\rho_2))(h) 
&= \frac {({}^{x_2}(\zeta_{x_1}(\rho_1)) \, \zeta_{x_2}(\rho_2))(h)(1)}
{\eta_{x_1x_2}(h)(1)} \\
&= \frac {(\zeta_{x_1}(\rho_1)(h))(x_2) \,\,\, (\zeta_{x_2}(\rho_2)(h))(1)}
{\eta_{x_1x_2}(h)(1)} \\
&= \frac {(\varphi_1(\rho_1)(h))(x_2) \,\, \eta_{x_1}(h)(x_2) \,\,
(\varphi_1(\rho_2)(h))(1) \,\, \eta_{x_2}(h)(1)} {\eta_{x_1x_2}(h)(1)} \\
&= (\nu(x_1, \, x_2) \, \rho_1^{x_2} \, \rho_2)(h) \\
\end{split}
\end{equation*}

\noindent Hence, $\overline{\zeta_{x_1x_2}(\nu(x_1, \, x_2) \, 
\rho_1^{x_2} \rho_2)} 
= \overline{{}^{x_2}(\zeta_{x_1}(\rho_1)) \, \zeta_{x_2}(\rho_2)}$.
This shows that T is a homomorphism. It should be clear that $T$ is an isomorphism 
and the theorem is proved.

\end{proof}

\begin{example}
Let $H$ be abelian and normal in $G$ such that its order is relatively prime 
to the order of the group $K$ and suppose $\psi({}^y\mu/\mu)$ is  
trivial in $H^2(H, \, k^\times)$, for all $y \in K$.
Then the Grothendieck ring of $Vec(G, \omega)^*_{\M(H, \, \mu)}$
is $\mathbb{Z}[\widehat{H} \rtimes K]$. 
Indeed, since $|H|$ and $|K|$ are relatively prime we have 
$H^2(K, \, \widehat{H}) = \{1\}$ which implies that 
$\nu$ is trivial in $H^2(K, \, \widehat{H})$. 

\end{example}

\end{subsection}


\begin{subsection}
{Skeleton of the category ${\C}_{\M}^{*}$}

A {\em skeleton} of a category $\D$ is any full subcategory $\overline{\D}$
such that each object of $\D$ is isomorphic (in $\D$) to exactly one
object of $\overline{\D}$. Every category is equivalent to any of its 
skeletons.
Let us recall how one constructs a skeleton $\overline{\D}$ 
of any tensor category $\D$ with associativity contraint $a$ and tensor
product $\otimes$. 
The construction is as follows: choose one object from each isomorphism
class of objects in $\D$.
Let $obj(\overline{\D})$ be the set of all objects choosen
above. For any $X \in obj(\D)$, by $\overline{X}$ we mean the object in
$\overline{\D}$ that represents the object $X$.

Define $Hom_{\overline{\D}}(X, \, Y)
= Hom_{\D}(X, \, Y)$. Define tensor product
$\odot$ in $\overline{\D}$: $X \odot Y
= \overline{X \otimes Y}$ for $X$, $Y \in Obj(\overline{\D})$.
Fix isomorphisms $\sigma(X, \, Y):
X \odot Y \tilde{\to} X \otimes Y$ in $\D$, for all $X$, $Y \in \overline{\D}$.
For any $f \in Hom_{\overline{\D}}(X, \, Y)$
and $g \in Hom_{\overline{\D}}(X^\prime, \, Y^\prime)$
define its tensor product:$f \odot g = \sigma(X^\prime, \, Y^\prime)^{-1}
\circ (f \otimes g) \circ \sigma(X, \, Y)$.

We now define associativity constraint $\overline{a}$ in
$\overline{\D}$. For any $X, \, Y, \,
Z \in obj(\overline{\D})$ define $\overline{a}_{X, \, Y, \,
Z}$ to be the following composition.  

\begin{equation*}
(X \odot Y) \odot Z
\xrightarrow{\sigma(X \odot Y, \, Z)}
(X \odot Y) \otimes Z
\xrightarrow{\sigma(X, \,  Y) \otimes id_{Z}}
(X \otimes Y) \otimes Z
\xrightarrow{a_{X, \, Y, \,Z}} 
X \otimes (Y \otimes Z)
\end{equation*}

\begin{equation*}
\xrightarrow{(id_{X} \otimes\sigma(Y, \,  Z))^{-1}} 
X \otimes (Y \odot Z)
\xrightarrow{\sigma(X, \,  Y \odot Z)^{-1}}
X \odot (Y \odot Z)
\end{equation*}

Left and right unit constraints are defined in the obvious way. It can be
shown that the necessary axioms (pentagon, triangle) are satisfied. Hence
$\overline{\D}$ is a monoidal category. One can also show that the categories
$\D$ and $\overline{\D}$ are tensor equivalent.

\begin{remark}
\label{ind}
If $\D$ is a pointed fusion category, then the simple objects of
$\overline{\D}$ form a group and the associativity constraint
in $\overline{\D}$ gives rise to a 3-cocycle. The cohomology class of 
this 3-cocycle does not depend on 
the choices made in the construction of $\overline{\D}$. 
\end{remark}

The function $\kappa$ defines an element in $Z^2(K, \, H)$:

\begin{equation}
\label{kappa1}
\kappa(x_1, \, x_2) := \kappa_{x_1, \, u(x_2)}. 
\end{equation}

\noindent Note that
the cohomology class of $\kappa$ is independent on the choice of
the function $u$. Also note that the cohomology class that $\kappa$ defines in
$H^2(K, \, H)$ is equal to the cohomology class associated to the 
the exact sequence $1 \to H \to G \to K \to 1$.

Define a 3-cocycle on the group $\widehat{H} \rtimes_{\nu} K$ 
with coefficients in $k^\times$
(see (\ref{cproduct1}) and (\ref{nutilde})):

\begin{equation}
\label{varpi}
\varpi((\rho_1, \, x_1), \, (\rho_2, \, x_2), \, (\rho_3, \, x_3)) =
(\tilde{\nu}(x_1, \, x_2)(u(x_3)))(1) \, \, \rho_1(\kappa(x_2, \, x_3))
\end{equation}

\noindent for any $(\rho_1, \, x_1)$, $(\rho_2, \, x_2)$, $(\rho_3, \, x_3)
\in \widehat{H} \rtimes_{\nu} K$. 

\begin{remark}
(i) It is routine to check that $\varpi$ does indeed define a 3-cocycle and
that its cohomology class does not depend on the choice of the 
function $u: K \to G$.

(ii) A special case, with $\tilde{\nu} \equiv 1$, of the formula in (\ref{varpi})
appeared in \cite{GMN}.
\end{remark}

\begin{theorem}
\label{thm3}
The fusion categories ${\C}_{\M}^{*}$ and $Vec(\widehat{H} \rtimes_{\nu} K,
\, \varpi)$ are equivalent.
\end{theorem}

\begin{proof}
Let us construct a skeleton $\overline{{\C}_{\M}^{*}}$ of the
category ${\C}_{\M}^{*}$.   
Let $\overline{\Lambda} = \bigcup_{x \in K} 
\{(\zeta_x(\rho), \, x) \, | \, \rho \in \widehat{H}\}$ 
denote the set of all simple objects of 
$\overline{{\C}_{\M}^{*}}$. See (\ref{zeta}) for definition of $\zeta_x$.
Tensor product $\odot$ in $\overline{{\C}_{\M}^{*}}$:
$(\zeta_{x_1}(\rho_1), \, x_1) \odot (\zeta_{x_2}(\rho_2), \, x_2)
= \overline{(\zeta_{x_1}(\rho_1), \, x_1) \otimes
(\zeta_{x_2}(\rho_2), \, x_2)} = \overline{({}^{x_2}(\zeta_{x_1}(\rho_1)) 
\, \zeta_{x_2}(\rho_2), \, x_1x_2)}
= (\zeta_{x_1x_2} (\nu(x_1, \, x_2) \rho_1^{x_2} \rho_2), \, x_1x_2)$. 
Note that $\overline{\Lambda}$ forms a group (multiplication coming from $\odot$) 
that is isomorphic to $\widehat{H} \rtimes_{\nu} K$.

Fix isomorphisms in ${\C}_{\M}^{*}$: $C \ni
f((\zeta_{x_1}(\rho_1), \, x_1), \, (\zeta_{x_2}(\rho_2), \, x_2)) :
(\zeta_{x_1}(\rho_1), \, x_1) \odot (\zeta_{x_2}(\rho_2), \, x_2)
= (\zeta_{x_1x_2}(\nu(x_1, \, x_2) \rho_1^{x_2} \rho_2), \, x_1x_2) \tilde{\to}
({}^{x_2}(\zeta_{x_1}(\rho_1)) \, \zeta_{x_2}(\rho_2), \, x_1x_2)
= (\zeta_{x_1}(\rho_1), \, x_1) \otimes (\zeta_{x_2}(\rho_2), \, x_2)$, for
all $(\zeta_{x_1}(\rho_1), \, x_1)$, $(\zeta_{x_2}(\rho_2), \, x_2) \in 
\overline{\Lambda}$. The following equality must hold: \\

\begin{equation*}
\begin{split}
{}^{x_2}(\zeta_{x_1}(\rho_1)) \, \zeta_{x_2}(\rho_2)) (g)(y) =
& \frac {f((\zeta_{x_1}(\rho_1), \, x_1), \, (\zeta_{x_2}(\rho_2), \, x_2))
(y \triangleleft g)} {f((\zeta_{x_1}(\rho_1), \, x_1), \, 
(\zeta_{x_2}(\rho_2), \, x_2))(y)} \\ \\
& \qquad \times {\zeta_{x_1x_2} (\nu(x_1, \, x_2) \rho_1^{x_2} \rho_2)(g)(y)} \\ \\
\end{split}
\end{equation*}

\noindent for all $g \in G$, $y \in K$. After using the definition of $\zeta_{x_1}$,
$\zeta_{x_2}$, $\zeta_{x_1x_2}$ and $\tilde{\nu}$, canceling and rearranging, 
the above relation becomes: \\

\begin{equation*}
\frac {f((\zeta_{x_1}(\rho_1), \, x_1), \, (\zeta_{x_2}(\rho_2), \, x_2))
(y \triangleleft g)} {f((\zeta_{x_1}(\rho_1), \, x_1), \, (\zeta_{x_2}(\rho_2), \, x_2))(y)}
= \frac {\tilde{\nu}(x_1, \, x_2)(g)(y) \, \, \rho_1(\kappa_{x_2y, \, g})} 
{\nu(x_1, \, x_2)(\kappa_{y, \, g}) \, \, \rho_1^{x_2}(\kappa_{y, \, g})} 
\end{equation*}

\vspace{0.1in}

\noindent Putting $y = 1$ and $g = u(y)$ in the above relation and canceling, 
we obtain: 

\begin{equation}
\label{a}
f((\zeta_{x_1}(\rho_1), \, x_1), \, (\zeta_{x_2}(\rho_2), \, x_2))(y)
=\frac {\tilde{\nu}(x_1, \, x_2)(u(y))(1) \, \, \rho_1(\kappa_{x_2, \, u(y)})}
{f((\zeta_{x_1}(\rho_1), \, x_1), \, (\zeta_{x_2}(\rho_2), \, x_2))(1)} 
\end{equation}

\vspace{0.1in}

Now let us calculate the associativity constraint in 
$\overline{{\C}_{\M}^{*}}$ which we denote
by $\varpi^\prime$. For any $(\zeta_{x_1}(\rho_1), \, x_1)$,  
$(\zeta_{x_2}(\rho_2), \, x_2)$, $(\zeta_{x_3}(\rho_3), \, x_3) 
\in \overline{\Lambda}$, $\varpi^\prime$ is defined as: \\

$\varpi^\prime((\zeta_{x_1}(\rho_1), \, x_1), \, (\zeta_{x_2}(\rho_2), \, x_2), \,
(\zeta_{x_3}(\rho_3), \, x_3))$

\begin{equation*}
\begin{split}
&= \frac {(f((\zeta_{x_1}(\rho_1), \, x_1), \, (\zeta_{x_2}(\rho_2), \, x_2)))
\otimes Id_{(\zeta_{x_3}(\rho_3), \, x_3)})}
{f((\zeta_{x_1}(\rho_1), \, x_1), \, (\zeta_{x_2}(\rho_2), \, x_2) \odot 
(\zeta_{x_3}(\rho_3), \, x_3))}  \\
& \qquad \qquad \qquad \qquad \times \frac {f((\zeta_{x_1}(\rho_1), \, x_1) 
\odot (\zeta_{x_2}(\rho_2), \, x_2),
\, (\zeta_{x_3}(\rho_3), \, x_3))}
{(Id_{(\zeta_{x_1}(\rho_1), \, x_1)} \otimes
f((\zeta_{x_2}(\rho_2), \, x_2), \, (\zeta_{x_3}(\rho_3), \, x_3)))} \\ \\
&= \frac {{}^{x_3}(f((\zeta_{x_1}(\rho_1), \, x_1), \, (\zeta_{x_2}(\rho_2), \, x_2)))}
{f((\zeta_{x_1}(\rho_1), \, x_1), \, (\zeta_{x_2}(\rho_2), \, x_2) \odot 
(\zeta_{x_3}(\rho_3), \, x_3))}  \\
& \qquad \qquad \qquad \qquad \times \frac {f((\zeta_{x_1}(\rho_1), \, x_1) \odot 
(\zeta_{x_2}(\rho_2), \, x_2),
\, (\zeta_{x_3}(\rho_3), \, x_3))}
{f((\zeta_{x_2}(\rho_2), \, x_2), \, (\zeta_{x_3}(\rho_3), \, x_3))} \\ \\
\end{split}
\end{equation*}

Note that $\varpi^\prime((\zeta_{x_1}(\rho_1), \, x_1), \, (\zeta_{x_2}(\rho_2), \, x_2), \, 
(\zeta_{x_3}(\rho_3), \, x_3))$ is an automorphism of
a simple object in ${\C}_{\M}^{*}$. By Note \ref{auto}, 
$\varpi^\prime((\zeta_{x_1}(\rho_1), \, x_1), \, 
(\zeta_{x_2}(\rho_2), \, x_2), \, (\zeta_{x_3}(\rho_3), \, x_3))(y)$ is constant for 
all $y \in K$. Thus, it suffices to calculate
$\varpi^\prime((\zeta_{x_1}(\rho_1), \, x_1), \, (\zeta_{x_2}(\rho_2), \, x_2), \,
(\zeta_{x_3}(\rho_3), \, x_3))(1)$. We have, \\ 

$\varpi^\prime((\zeta_{x_1}(\rho_1), \, x_1), \, (\zeta_{x_2}(\rho_2), \, x_2), \,
(\zeta_{x_3}(\rho_3), \, x_3))(1)$

\begin{equation*}
\begin{split}
&= \frac {f((\zeta_{x_1}(\rho_1), \, x_1), \, (\zeta_{x_2}(\rho_2), \, x_2))(x_3)}
{f((\zeta_{x_1}(\rho_1), \, x_1), \, (\zeta_{x_2}(\rho_2), \, x_2) \odot (\zeta_{x_3}(\rho_3), \, x_3))(1)}  \\
&\qquad \qquad \qquad \qquad \times \frac
{f((\zeta_{x_1}(\rho_1), \, x_1) \odot (\zeta_{x_2}(\rho_2), \, x_2),
\, (\zeta_{x_3}(\rho_3), \, x_3))(1)}
{f((\zeta_{x_2}(\rho_2), \, x_2), \, (\zeta_{x_3}(\rho_3), \, x_3))(1)} \\
&= \frac {f((\zeta_{x_1}(\rho_1), \, x_1), \, (\zeta_{x_2}(\rho_2), \, x_2))(1)}
{f((\zeta_{x_1}(\rho_1), \, x_1), \, (\zeta_{x_2}(\rho_2), \, x_2) \odot (\zeta_{x_3}(\rho_3), \, x_3))(1)}  \\
&\qquad \qquad \qquad \qquad \times \frac
{f((\zeta_{x_1}(\rho_1), \, x_1) \odot (\zeta_{x_2}(\rho_2), \, x_2),
\, (\zeta_{x_3}(\rho_3), \, x_3))(1)}
{f((\zeta_{x_2}(\rho_2), \, x_2), \, (\zeta_{x_3}(\rho_3))), \, x_3)(1)} \\
&\qquad \qquad \qquad \qquad \times\tilde{\nu}(x_1, \, x_2)(u(x_3))(1) \, \,
\rho_1(\kappa_{x_2, \, u(x_3)}) \\ \\
\end{split}
\end{equation*}

\noindent We used (\ref{a}) to obtain the last equality.

Since the cohomology class of $\varpi^\prime$ does not depend on the choice of
the isomorphisms $f(\cdot, \, \cdot)$, we can assume that $f(\cdot, \, \cdot)(1) = 1$. 
Also, regard $\varpi^\prime$ as a 3-cocycle on $\widehat{H} \rtimes_{\nu} K$. Then we get:

\begin{equation*}
\varpi^\prime((\rho_1, \, x_1), \, (\rho_2, \, x_2), \, (\rho_3, \, x_3)) =
\tilde{\nu}(x_1, \, x_2)(u(x_3))(1) \, \, \rho_1(\kappa(x_2, \, x_3))
\end{equation*}

\noindent for any $(\rho_1, \, x_1)$, $(\rho_2, \, x_2)$, $(\rho_3, \, x_3)
\in \widehat{H} \rtimes_{\nu} K$. That is, $\varpi^\prime = \varpi$ and the theorem
is proved.
\end{proof}

\begin{example}
Let $G = \mathbb{Z}/4\mathbb{Z} = 
\{\overline{0}, \overline{1}, \overline{2}, \overline{3}\}$, $\omega = 1$, 
$H =\{\overline{0}, \overline{2}\}$, and $\mu \equiv 1$.
Since $\mu \equiv 1$ we can assume that $\tilde{\nu} \equiv 1$ 
(see (\ref{nutilde})) and $\nu \equiv 1$ (see Corollary \ref{asso}).
By Theorem \ref{thm2} it follows that $\mathcal{K}_0({\C}_{\M}^{*}) \cong 
\mathbb{Z}[\widehat{\mathbb{Z}/2\mathbb{Z}} \times 
\mathbb{Z}/2\mathbb{Z}]$. Let $\widehat{\mathbb{Z}/2\mathbb{Z}} = 
\{\rho_0, \, \rho_1\}$,
where $\rho_1$ represents the non-trivial character. We have, 
$K = \{H+\overline{0}, \, H+\overline{1}\}$. We claim that the 
associativity constraint $\varpi$ in  
${\C}_{\M}^{*}$ is non-trivial. It suffices to show that the
restriction of $\varpi$ to some non-trivial subgroup of 
$\widehat{\mathbb{Z}/2\mathbb{Z}} \times 
\mathbb{Z}/2\mathbb{Z}$ is non-trivial. Consider the 
restriction of $\varpi$ to the subgroup $K = \{(\rho_0, \, H+\overline{0}), 
(\rho_1, \, H+\overline{1})\}$. It suffices to 
show that there exists a triple of elements in this subgroup such
that $\varpi$ evaluated at this triple is not equal to $1$. Define the
function $u: K \to G$ by $u(H+\overline{0}) = \overline{0}$
and $u(H+\overline{1}) = \overline{1}$. Since $\mu \equiv 1$, we can choose
$\tilde{\nu} \equiv 1$. So the first factor in the definition of $\varpi$
vanishes. We have, $\varpi((\rho_1, \, H+\overline{1}), \, 
(\rho_1, \, H+\overline{1}), \, (\rho_1, \, H+\overline{1})) 
= \rho_1(\kappa(H+\overline{1}, \, H+\overline{1}))
= \rho_1(\overline{2}) = -1$. Thus, the 3-cocycle $\varpi$ is non-trivial.
In particular, the fusion categories
$Vec(\mathbb{Z}/4\mathbb{Z}, \, 1)$ and 
$Vec(\mathbb{Z}/2\mathbb{Z} \times \mathbb{Z}/2\mathbb{Z}, \, \varpi)$ are
weakly Morita equivalent (see next section).
\end{example}

\end{subsection}


\end{section}


\begin{section}
{Categorical Morita Equivalence}

Two tensor categories $\C$ and $\D$ are said to be {\em weakly 
Morita equivalent} if there exists an indecomposable right module 
category $\M$ over
$\C$ such that the categories ${\C}_{\M}^{*}$ and $\D$ 
are tensor equivalent (see \cite{Mu}). It was shown by M\"uger that this
indeed is an equivalence relation.

Using the notion of weak Morita equivalence 
we put an equivalence relation on the set of all pairs $(G, \, \omega)$, 
where $G$ is a finite group and $\omega \in H^3(G, \, k^\times)$:

\begin{definition}
We say that two pairs $(G, \, \overline{\omega})$ and $(G^\prime, \, 
\overline{\omega}^\prime)$ 
are {\em categorically Morita equivalent}
and write $(G, \, \overline{\omega}) \morita (G^\prime, \, \overline{\omega}^\prime)$
if the tensor categories $Vec(G, \, \omega)$ and
$Vec(G^\prime, \, \omega^\prime)$, are weakly Morita equivalent.
\end{definition}

\begin{remark}
Note that finding categorically Morita equivalence classes of the set 
of all pairs $(G, \, \omega)$, where $G$ is a finite
group and $\omega \in H^3(G, \, \, k^\times)$ amounts to finding 
weakly Morita equivalence classes of the set of all group-theoretical categories.
\end{remark}

\noindent We also define an equivalence relation on the set of all groups:

\begin{definition}
We say that two groups $G$ and $G^\prime$ are {\em categorically Morita equivalent}
and write $G \morita G^\prime$ if the pairs $(G, \, 1)$ and
$(G^\prime, \, 1)$ are categorically Morita equivalent.
\end{definition} 

\begin{remark}
Two finite groups $G$ and $G^\prime$ are called {\em isocategorical} if
their representation categories Rep($G$) and Rep($G^\prime$) 
are tensor equivalent \cite{EG}.
If two groups $G$ and $G^\prime$ are isocategorical, then they are
categorically Morita equivalent (this follows from the fact that
for any group $G$ the categories Rep($G$) and $Vec(G, \, 1)_{\M(G, \, 1)}^*$
are tensor equivalent). We show in Section \ref{ex} that the
converse is not true, that is, there do exist
groups that are categorically Morita equivalent but not isocategorical.
\end{remark} 

\begin{remark}
\label{center}
It was shown in \cite{O2} that if the tensor categories $\C$ and $\D$ are weakly 
Morita equivalent, then their Drinfeld centers are equivalent as braided tensor
categories. It follows that if two
groups are categorically Morita equivalent, then the Drinfeld centers of their 
representation categories are equivalent as braided tensor categories.
\end{remark}

\begin{definition}
We say that a group $G$ is categorically Morita rigid if any group that
is categorically Morita equivalent to $G$ is actually isomorphic to $G$.
\end{definition}

\begin{remark}
By remark \ref{center} it follows that abelian groups are categorically Morita 
rigid. In particular, an abelian group can not be categorically Morita 
equivalent to a non-abelian group.
\end{remark}

\noindent The next theorem gives a group-theoretical and cohomological 
interpretation of categorical Morita equivalence.

\begin{theorem}
\label{thm4}
Two pairs $(G, \, \overline{\omega})$ and 
$(G^\prime, \, \overline{\omega}^\prime)$ 
are categorically Morita equivalent if and only if the following
conditions hold:
\begin{enumerate}
\item 
$G$ contains a normal abelian subgroup $H$ such that 
$\omega|_{H \times H \times H}$ is trivial in $H^3(H, \, k^\times)$,
\item 
there is a 2-cochain $\mu \in C^2(G, \, C)$ such that 
$\delta^2\mu = \omega$ and $\psi({}^y\mu/\mu)$ is trivial
in $H^2(H, \, k^\times)$, for all $y \in H \setminus G$
and there is an isomorphism a: $G^\prime \tilde{\to} \widehat{H} \rtimes_{\nu} 
(H \setminus G)$ (see (\ref{cproduct1})) and
\item 
the $3$-cocycle $\frac {\varpi^a} {\omega^\prime}$ is 
trivial in $H^3(G^\prime, \, k^\times)$ 
(see (\ref{varpi}) for definition of $\varpi$).
\end{enumerate}
\end{theorem}

\begin{proof}
Suppose the pairs $(G, \, \overline{\omega})$ and 
$(G^\prime, \, \overline{\omega}^\prime)$ 
are categorically Morita equivalent. 
Then there exists an indecomposable right module category $\M$ over 
$Vec(G, \, \omega)$ such that the categories $Vec(G, \, \omega)_{\M}^*$ and 
$Vec(G^\prime, \, \omega^\prime)$ are tensor equivalent. So there exists
a subgroup $H$ of $G$ such that 
$\omega|_{H \times H \times H}$ represents the trivial class in $H^3(H, \, K^\times)$
and 2-cochain $\mu \in C^2(G, \, C)$ (satisfying the  
$\delta^2\mu = \omega$) 
which together produce the module category $\M$. Note that 
$Vec(G, \, \omega)_{\M}^*$ must be pointed. 
By Theorem \ref{thm1}, it follows that $H$ is abelian and normal in $G$ and that
$\psi({}^y\mu/\mu)$ is trivial
in $H^2(H, \, k^\times)$, for all $y \in H \setminus G$. Theorem \ref{thm3} says that 
$Vec(G, \, \omega)_{\M}^* \cong Vec(\widehat{H} \rtimes_{\nu} 
(H \setminus G), \, \varpi)$.
It now follows that there must exist an isomorphism 
$a: G^\prime \to \widehat{H} \rtimes_{\nu} (H \setminus G)$ such that
$\varpi^a$ is cohomologous to $\omega^\prime$.
The converse should be clear and the theorem is proved.
\end{proof}

\begin{corollary}
\label{cor1}
Two groups $G$ and $G^\prime$ 
are categorically Morita equivalent if and only if the following conditions hold:

\begin{enumerate}
\item
$G$ contains a normal abelian subgroup $H$, 
\item
there exists a $G$-invariant $\mu \in H^2(H, \, k^\times)$ such that
the groups $G^\prime$ and $\widehat{H} \rtimes_{\nu} 
(H \setminus G)$ are isomorphic (where $\nu = \Phi(\mu)$, see (\ref{Phi})) and
\item
the $3$-cocycle $\varpi$ defined in (\ref{varpi}) is trivial.
\end{enumerate}
\end{corollary}

\end{section}


\begin{section}
{Examples of categorically Morita equivalent groups with
non-isomorphic Grothendieck rings}
\label{ex}
In this section we produce a series of pairs of groups that 
are categorically Morita equivalent but have non-isomorphic Grothendieck rings.
Let $p$ and $q$ be odd primes such
that $p-1$ is divisible by $q$. Then there exists a unique upto isomorphism
non-trivial semidirect
product of the groups $\mathbb{Z}/p\mathbb{Z}$ and $\mathbb{Z}/q\mathbb{Z}$.
Let $a$ and $b$ be generators of the groups 
$\mathbb{Z}/p\mathbb{Z}$
and $\mathbb{Z}/q\mathbb{Z}$, respectively. 
Let us fix an action of
$\mathbb{Z}/q\mathbb{Z}$ on $\mathbb{Z}/p\mathbb{Z}$:
fix a $t \in \mathbb{Z}$ ($t \not\equiv 1 \mod p$) such that $t^q-1$ is divisible by $p$.
Such a $t$ of course exists because $p-1$ is divisible by $q$.
Then the action of $\mathbb{Z}/q\mathbb{Z}$ on $\mathbb{Z}/p\mathbb{Z}$
is defined by: 
$a \triangleleft b := a^t$.
Let $\rho$ be a generator of the groups $\widehat{\mathbb{Z}/p\mathbb{Z}}$.
Then the induced action of
$\mathbb{Z}/q\mathbb{Z}$ on $\widehat{\mathbb{Z}/p\mathbb{Z}}$ is given by: 
$(\rho \triangleleft b)(a) := \rho(a \triangleleft b^{-1})$.
But $b^{-1} = b^{q-1}$. So $\rho \triangleleft b = \rho^{t^{q-1}}$.

The subgroup
$\mathbb{Z}/p\mathbb{Z}$ (identified with $\mathbb{Z}/p\mathbb{Z} \times \{1\}$) 
of $\mathbb{Z}/p\mathbb{Z} \rtimes \mathbb{Z}/q\mathbb{Z}$
can be considered as a right 
($\mathbb{Z}/p\mathbb{Z} \rtimes \mathbb{Z}/q\mathbb{Z}$)-module where the action
is via conjugation. The dual group $\widehat{\mathbb{Z}/p\mathbb{Z}}$ is also
a right ($\mathbb{Z}/p\mathbb{Z} \rtimes \mathbb{Z}/q\mathbb{Z}$)-module with the
action being induced from the action of 
$\mathbb{Z}/p\mathbb{Z} \rtimes \mathbb{Z}/q\mathbb{Z}$ on $\mathbb{Z}/p\mathbb{Z}$.
Let $G := \mathbb{Z}/p\mathbb{Z} \rtimes (\mathbb{Z}/p\mathbb{Z} 
\rtimes \mathbb{Z}/q\mathbb{Z})$ and
$G^\prime := \widehat{\mathbb{Z}/p\mathbb{Z}} \rtimes (\mathbb{Z}/p\mathbb{Z} 
\rtimes \mathbb{Z}/q\mathbb{Z})$.

\begin{lemma}
The groups $G$ and $G^\prime$ have different number of normal subgroups of
order $p$.
\end{lemma}

\begin{proof}
Note that both groups have the same number of subgroups of order $p$.
We claim that all subgroups of order $p$ in $G$ are normal whereas there
exists a non-normal subgroup of order $p$ in $G^\prime$.
The generator of any subgroup of $G$ of order $p$ is of the form
$(a^l, \, (a^m, \, 1))$, where $l$, $m \in \{1, \dots, p\}$ with
$l$ and $m$ not simultaneously equal to $p$. 
The elements $(a, \, (1, \, 1))$,
$(1, \, (a, \, 1))$, and $(1, \, (1, \, b))$ generate the group $G$.
Note that the element $(a^l, \, (a^m, \, 1))$ is stable under conjugation
by the first two generators of $G$. While conjugation
by the third generator gives: $(1, \, (1, \, b))^{-1}
(a^l, \, (a^m, \, 1))(1, \, (1, \, b)) = (1, \, (1, \, b^{-1}))
(a^{lt}, \, (a^{mt}, \, b)) = (a^{lt}, \, (a^{mt}, \, 1)) = 
(a^l, \, (a^m, \, 1))^t$.
This shows that all subgroups of order $p$ in $G$ are normal. Consider the subgroup
of $G^\prime$ of order $p$ generated by the element $(\rho, \, (a, \, 1))$. 
We have \\
$(1, \, (1, \, b))^{-1}(\rho, \, (a, \, 1))(1, \, (1, \, b))
=(1, \, (1, \, b^{-1}))(\rho^{t^{q-1}}, \, (a^t, \, b)) = 
(\rho^{t^{q-1}}, \, (a^t, \, 1))$. 
Note that the element $(\rho^{t^{q-1}}, \, (a^t, \, 1))$ is not a power of
$(\rho, \, (a, \, 1))$ because $t^{q-1} \not\equiv t \mod p$.
This shows that the subgroup of $G^\prime$ of order $p$ 
generated by the element $(\rho, \, (a, \, 1))$ is not normal and hence 
the lemma is proved.
\end{proof}

\begin{corollary}
\label{thm5}
The groups $G$ and $G^\prime$ are categorically 
Morita equivalent but have non-isomorphic Grothendieck rings.
\end{corollary}

\begin{proof}
To see that these two groups satisfy the conditions
in Corollary \ref{cor1}, take $H$ to be the subgroup $\mathbb{Z}/p\mathbb{Z}$
of $G$ and take $\mu \equiv 1$. Observe that 
the groups $\widehat{H} \rtimes (H \setminus G)$ and $G^\prime$ are isomorphic.
Since the exact sequence 
$1 \to H \to G \to H \setminus G \to 1$ splits we can assume $\kappa \equiv 1$.
Also, since $\mu \equiv 1$ we can assume that $\tilde{\nu} \equiv 1$
and therefore $\varpi \equiv 1$.
It follows that the groups 
$G$ and $G^\prime$ are categorically Morita equivalent.
To see that these groups have non-isomorphic Grothendieck rings
note that the above lemma implies that these groups have different number of 
quotient groups of order $pq$. By \cite[Proposition 3.11]{N} it follows that the
Grothendieck rings $K_0$(Rep($G$)) and $K_0$(Rep($G^\prime$)) of the 
two groups are not isomorphic.
\end{proof}

\begin{corollary}
The groups $G$ and $G^\prime$ are non-isocategorical.
\end{corollary}

\begin{proof}
This follows immediately from the above corollary.
\end{proof}

\begin{remark}
(i) By Remark \ref{center} the representation categories Rep($D(G)$) and 
Rep($D(G^\prime)$) of the Drinfeld doubles of the groups $G$ and $G^\prime$
are equivalent as braided tensor categories and hence these groups define the same 
modular data. Note also that the Hopf algebras $D(G)$ and $D(G^\prime)$ must be
guage-equivalent. \\

(ii) Equivalence of certain twisted doubles of groups was investigated in \\
\cite{GMN}. \\

(iii) The above examples of categorically Morita equivalent
groups come from a more general construction: start with any finite group
$G$ and a finite right $G$-module $H$. Consider the semidirect product $H \rtimes G$.
We can regard $\widehat{H}$ as a right 
$G$-module with the action being induced from the action of
$G$ on H. Then the groups $H \rtimes G$ and
$\widehat{H} \rtimes G$ are categorically Morita equivalent. Note
however that these two groups are not always non-isomorphic. \\

(iv) By Ito's theorem \cite[Theorem 6.3.9]{G} it follows
that the possible dimensions of irreducible representations of the groups
$G$ and $G^\prime$
are $1$ and $q$. It can be shown that the order of the commutator subgroup is $p^2$
for both groups. Therefore, the order of the abelianization (equal to the number
of $1$-dimensional representations) of both groups is $q$. So the group algebras $kG$ 
and $kG^\prime$ are both isomorphic to 
$\underbrace{k \oplus k \oplus \dots \oplus k}_{q \,\, copies} \oplus 
\underbrace{M_q(k) \oplus M_q(k) \oplus \dots \oplus 
M_q(k)}_{(p^2-1)/q \,\, copies}$. \\

(v) It follows from Corollary \ref{thm5} that the groups $G$ and $G^\prime$ 
have different character tables. This provides
a counter-example to the hunch mentioned in \cite{CGR} that groups defining the 
same modular data will have the same character table.

\end{remark}

\end{section}


\setcounter{secnumdepth}{0}\begin{section}
{Acknowledgments}
The author would like to express his deep gratitude to his dissertation advisor, 
Prof. Dmitri Nikshych, for introducing this problem to him and for providing 
invaluable help and guidance during the preparation of this paper. The author
would also like to thank the referee for comments that helped improve the paper.
Thanks are also due to Shamindra Ghosh for useful discussions. 
During preparation of the paper the author was supported by fellowships awarded by
the graduate school at the University of New Hampshire and by NSF grant 
DMS-0200202.
\end{section}

\bibliographystyle{ams-alpha}


\end{document}